\documentclass[graybox]{svmult}

\usepackage{type1cm}        % activate if the above 3 fonts are
\usepackage{graphicx}        % standard LaTeX graphics tool
\usepackage{multicol}        % used for the two-column index
\usepackage[bottom]{footmisc}% places footnotes at page bottom

\usepackage{newtxtext}       %
\usepackage{newtxmath}       % selects Times Roman as basic font

\usepackage[english]{babel}
\usepackage[babel=true]{microtype}
\usepackage{amsmath,bm}

\usepackage{siunitx}
\usepackage{hyperref}
\usepackage{cleveref}
\crefname{equation}{}{}
\Crefname{equation}{Equation}{Equations}
\crefname{figure}{Fig.}{Figs.}
\Crefname{figure}{Figure}{Figures}
\crefname{table}{Table}{Tables}
\crefname{algorithm}{Algorithm}{Algorithms}
\crefname{section}{Sect.}{Sects.}
\Crefname{section}{Section}{Sectsions}

\crefname{theorem}{Theorem}{Theorems}
\Crefname{theorem}{Theorem}{Theorems}
\crefname{corollary}{Corollary}{Corollaries}
\Crefname{corollary}{Corollary}{Corollaries}
\crefname{lemma}{Lemma}{Lemmas}
\Crefname{lemma}{Lemma}{Lemmas}
\crefname{definition}{Definition}{Definitions}
\crefname{definition}{Definition}{Definitions}
\crefname{remark}{Remark}{Remarks}
\Crefname{remark}{Remark}{Remarks}

\newcommand{\Adiff}{A_{\mathrm{diff}}}
\newcommand{\Cc}{\mathbb{C}}
\newcommand{\gbc}{g_{\mathrm{bc}}}
\newcommand{\geqs}{\geqslant}
\newcommand{\gpeak}{g_{\mathrm{peak}}}
\newcommand{\leqs}{\leqslant}
\newcommand{\Rr}{\mathbb{R}}
\newcommand{\Rrnn}{\Rr^{N\times N}}
\newcommand{\tol}{\texttt{tol}}
\newcommand{\yex}{y_{\mathrm{ex}}}
\newcommand{\yref}{y_{\mathrm{ref}}}

\makeindex             % used for the subject index
\begin{document}

\title*{Exponential time integrators for unsteady advection--diffusion problems on refined meshes}
\titlerunning{Exponential time integrators for unsteady convection--diffusion}
\author{M.A. Botchev}
\institute{%
  M.A. Botchev \at{}
  Keldysh Institute of Applied Mathematics, 
  Russian Academy of Sciences,    
  Miusskaya Sq.~4, Moscow 125047, Russia,
and Marchuk Institute of Numerical Mathematics, Russian Academy of Science,
Gubkina St.~8, 119333 Moscow, Russia.
This work is supported by the Russian Science
Foundation grant No.~19-11-00338,
  \email{botchev@ya.ru}
   %\and{}
   %Name of Second Author \at{} Name, Address of Institute \email{name@email.address}
}
\maketitle

\abstract*{Time integration of advection dominated
advection--diffusion problems
on refined meshes can be a challenging task, since local refinement
can lead to a severe time step restriction, whereas standard
implicit time stepping is usually hardly suitable for treating
advection terms. We show that exponential time integrators
can be an efficient, yet conceptually simple, option in this case.
Our comparison includes three exponential integrators
and one conventional scheme, the two-stage Rosenbrock method ROS2
which has been a popular alternative to splitting methods
for solving advection--diffusion problems.} 

\abstract{Time integration of advection dominated
advection--diffusion problems
on refined meshes can be a challenging task, since local refinement
can lead to a severe time step restriction, whereas standard
implicit time stepping is usually hardly suitable for treating
advection terms. We show that exponential time integrators
can be an efficient, yet conceptually simple, option in this case.
Our comparison includes three exponential integrators
and one conventional scheme, the two-stage Rosenbrock method ROS2
which has been a popular alternative to splitting methods
for solving advection--diffusion problems.}

%---------------------------------------------------------------------
\section{Introduction}\label{sec:1}
Time integration of unsteady advection--diffusion
problems discretized in space on locally refined meshes
can be a challenging problem.
This is especially the case for advection dominated problems.
On the one hand, requirements of accuracy, monotonicity and total variation
diminishing (TVD) usually rule out the use of implicit time
integration for advection terms~\cite[Chapter~III.1.3]{HundsdorferVerwer:book}.
On the other hand, locally refined meshes can impose a severe
CFL stability restriction on the time step, thus making
explicit schemes very inefficient.

Within the method of lines framework, i.e., when discretization
in space is followed by time integration, different approaches exist
to cope with this problem.
A straightforward and widely used approach is
operator splitting~\cite[Chapter~IV]{HundsdorferVerwer:book},\cite[Chapter~3]{Zlatev},
when advection is usually treated explicitly in time and diffusion implicitly.
Though being conceptually simple and easy to apply in practice,
splitting methods unavoidably lead to splitting
errors, see, e.g.,~\cite{LanserVerwer1999,Csomos_ea2005}.
Moreover, a proper use of boundary conditions within
the splitting is sometimes not trivial and may lead to
error order reduction~\cite{Sommeijer_ea1981,EinkemmerOstermann2015}.
To reduce splitting errors many various approaches have
been proposed, with none being fully successful.
Here we mention source splitting techniques~\cite{jansurvey}
and Rosenbrock %AMF (Approximate Matrix Factorization)
schemes~\cite[Chapter~IV.5]{HundsdorferVerwer:book}.

Other possible approaches to integrate advection--diffusion problems
efficiently include implicit-explicit (IMEX)
methods~\cite[Chapter~IV.4]{HundsdorferVerwer:book} and
multirate
schemes~\cite{Savcenco_ea2007,ConstantinescuSandu2007,Schlegel_ea2009}.
% The latter are conceptually more complicated clearly require serious
% Splitting schemes
% multirate schemes 

In this paper we show that in some cases exponential time integration
schemes can serve as an efficient yet simple way to
integrate advection--diffusion problems in time on locally refined meshes.
Similar to implicit schemes, exponential schemes have attractive
stability properties.  However, exponential schemes have also excellent
accuracy properties and in some cases, especially for linear
ODE (ordinary differential equation) systems, are able to produce exact
solution to initial-value problem (IVP) being solved.  This is the
property we exploit in this work.
An example of an IVP that can be solved exactly by an exponential
solver is
\begin{equation}
\label{ivp0}
y'(t) = -Ay(t) + g, \quad y(0)=v,\quad t\in[0,T],
\end{equation}
where $v,g\in\Rr^N$ are given and $A\in\Rrnn$ represents the
advection--diffusion operator discretized in space.
Exponential solvers involve matrix-vector products with the matrix exponential
and related, the so-called $\varphi$ functions.  More specifically,
the exact solution $y(t)$ of~\eqref{ivp0} can be written as
\begin{equation}
\label{yex}
y(t) = v + t\varphi(-tA)(g-Av), \quad t\geqs 0,  
\end{equation}
where a matrix-vector product of the matrix function
$\varphi(-tA)$ and the vector $g-Av$ has to be computed.
Here the function $\varphi$ is defined as~\cite{HochbruckOstermann2010}
$$
\varphi(z)=
\begin{cases}
  \dfrac{e^z-1}{z}, \quad &\text{for } z\ne0, \; z\in\Cc,\\
  1,                      &\text{for }z=0.
\end{cases}
$$
In case $g\equiv 0$ expression~\eqref{yex} reduces to a familiar relation
$$
y(t)= \exp(-tA)v, \quad t\geqs 0.  
$$ 
Note that such an ``all-at-once'' exact solution of
systems~\eqref{ivp0} is also possible if $g$ is a given vector
function $g(t)$ of time $t$, see~\cite{Botchev2013}. %???
To solve general nonlinear IVPs, across-time iterations of the
waveform relaxation type can successfully be employed~\cite{Kooij_ea2017}.
Then, at each waveform relaxation iteration, a linear IVP of the
type~\eqref{ivp0}
is solved by an exponential scheme.  Such an approach is attractive
because no time stepping is involved and solution can often be obtained
for the whole time interval $t\in[0,T]$.

In this work we present comparison results for three
different exponential solvers.
All these methods are based on the Krylov subspace techniques
discussed in
\cite{EXPOKIT,Botchev2013,ART}.
Krylov subspace methods have been successfully used for evaluating
matrix exponential and related matrix functions
since the eighties, see in chronological
order~%
\cite{ParkLight86,Henk:f(A),DruskinKnizh89,Knizh91,Saad92,DruskinKnizh95,HochLub97}.
An attractive property of the Krylov subspace methods, which distinguishes them
from the other methods used for large matrix function evaluations $f(A)v$,
is their adaptivity with respect to the spectral properties of $A$
and the particular vector~$v$, see~\cite{Henk:book}.
% Eventually say about other than Krylov methods for expm
To work efficiently, Krylov subspace methods often need
a restarting~\cite{EiermannErnst06,GuettelFrommerSchweitzer2014},
a mechanism allowing to keep Krylov subspace dimension restricted
while preserving convergence of the unrestarted method.
% exponential solvers based on Krylov subspaces.   

The structure of this paper is as follows.
In Section~\ref{s:methods} the problem and methods used for its
solution are presented.
Section~\ref{s:num} is devoted to numerical experiments,
here the methods are compared and comparison results are
discussed.  Finally, some conclusions are drawn in Section~\ref{s:concl}.

%---------------------------------------------------------------------
\section{Problem formulation and methods}
\label{s:methods}
In this paper we assume that a linear PDE of the advection--diffusion type is
solved by the method of lines and, after a suitable space discretization,
the following IVP has to be solved
\begin{equation}
\label{ivp}
y'(t) = -Ay(t) + g(t), \quad y(0)=v,\quad t\in[0,T].
\end{equation}
Here $A\in\Rrnn$ represents the discretized advection--diffusion
operator and the given vector function $g(t)$ accounts for
time dependent sources or boundary conditions.

\subsection{Exponential time integrators}
Perhaps the simplest exponential integrator is exponential Euler
method which, applied to problem~\eqref{ivp}, reads
\begin{equation}
y_{n+1}=y_n + \Delta t\varphi(-\Delta t A)(g_n-Ay_n),
\end{equation}
where $\Delta t$ is the time step size, $y_n$ is numerical
solution at time $t=\Delta t n$ and $g_n=g(\Delta t n)$.
The method is inspired by relation~\eqref{yex} and for constant source
term $g$ is exact.  It is not difficult to check that it is first
order accurate.

Using extrapolation~\cite{BotchevVerwer09,Zlatev_ea2018book},
i.e., by combining solutions obtained with
different time steps, higher order methods can be obtained.
In~\cite{Botchev2013} globally extrapolated second order exponential
Euler method (EE2) is considered.  Its combination with EXPOKIT~\cite{EXPOKIT},
used to evaluate the $\varphi$ matrix vector products, argued to
be a competitive integrator.
The \texttt{phiv} function of EXPOKIT is able to efficiently
compute actions of the $\varphi$ matrix functions by a restarted
Arnoldi process, where the restarting is done by time stepping based on an error
estimation.
In experiments presented below we show that EE2/EXPOKIT can be significantly
improved by introducing residual-based error
control~\cite{CelledoniMoret97,DruskinGreenbaumKnizhnerman98,BGH13}
and replacing the restarting procedure by the residual-time (RT)
restarting presented in~\cite{ART}.

EE2 is an exponential integrator which evaluates matrix functions
\emph{within} a time stepping procedure: at each time step
$\varphi$ is computed by a Krylov subspace method.  It is often more
efficient~\cite{Botchev2016,BotchevOseledetsTyrtyshnikov2014},
if possible, to organize work in such a way that numerical linear algebra
work for matrix function evaluations is done \emph{across} time
stepping, for a certain time interval.
For instance, as shown in~\cite{Botchev2013}, if
for a certain time range $g(t)$
allows an approximation
\begin{equation}
\label{g=Up}  
g(t)\approx Up(t), \quad U\in\Rr^{N\times m},\quad m\ll N,
\end{equation}
then~\eqref{ivp} can be solved for this time range by a single
projection on a block Krylov subspace.

To be more specific, consider problem~\eqref{ivp0}, where
for simplicity and without loss of generality assume $v=0$.
A usual Krylov subspace solution of~\eqref{ivp0} constructs
a matrix $V_k\in\Rr^{N\times k}$, whose orthonormal columns span
the Krylov subspace~\cite{Henk:book}
$$
\mathrm{span} (g, Ag, \dots, A^{k-1}g), 
$$
and, searching for an approximate solution $y_k(t)=V_ku(t)\approx y(t)$, 
reduces~\eqref{ivp0} to its Galerkin projection
\begin{equation}
\label{ivp_u}  
V_k^TV_ku'(t)= -V_k^TAV_ku(t) + V_k^Tg
\quad\Leftrightarrow\quad
u'(t) = -H_ku(t) + \beta e_1,
\end{equation}
where $H_k=V_k^TAV_k$, $e_1=(1,0,\dots,0)^T\in\Rr^k$ is the first canonical
basis vector and $\beta=\|g\|$.  We have $V_k^Tg=V_k^TV_k(\beta e_1)=\beta e_1$
because, by construction, the first column of $V_k$ is $g/\|g\|$.
The small projected IVP~\eqref{ivp_u} can be solved by relation~\eqref{yex},
evaluating $\varphi(-t H_k)$ with well developed matrix function
techniques for small matrices (see, e.g.,~\cite{Higham_bookFM}).

Now consider, still assuming $v=0$, problem~\eqref{ivp} where $g(t)$
allows~\eqref{g=Up}.  Projecting~\eqref{ivp} on a block
Krylov subspace~\cite{Henk:book}
$$
\mathrm{span} (U, AU, \dots, A^{k-1}U), 
$$
we can reduce~\eqref{ivp} to its projected form
$$
u'(t) = -H_{k}u(t) + E_1 p(t),
$$
where we now have $H_{k}\in\Rr^{km\times km}$ and
$E_1\in\Rr^{km\times m}$ is a matrix whose columns  
are the first $m$ columns of the $km\times km$ identity matrix.
These observations lead to the exponential block Krylov (EBK)
method described in~\cite{Botchev2013}.  Note that approximation~\eqref{g=Up}
can be easily constructed by truncated singular value decomposition (SVD)
of the vectors $g(t_i)$, at a small additional cost~\cite{Botchev2013}.
This procedure simultaneously provides an error estimation
in~\eqref{g=Up}, so that a proper value for $m$ can be chosen.

The EBK solver exploits a stopping criterion and restarting which are
based on the exponential residual
concept~\cite{CelledoniMoret97,DruskinGreenbaumKnizhnerman98,BGH13}.
In particular, EBK iterations stop as soon as for the computed approximate
solution $y_k(t)$ holds
$$
\|r_k(t)\|\leqs\tol,\quad
r_k(t)\equiv -Ay_k(t)+g(t) - y_k'(t),\quad t\in[0,T].
$$

\subsection{ROS2 method: beyond splitting}
Rosenbrock schemes~\cite[Chapter~IV.5]{HundsdorferVerwer:book}
have been a popular alternative to
splitting methods, as they allow to reduce splitting errors
and avoid other negative effects related to splitting, such as
order reduction.  Let $f(t,y)=-Ay(t)+g(t)$ be the ODE right hand side
in~\eqref{ivp}.  The two-stage Rosenbrock method ROS2 reads
\begin{equation}
\label{ros2}
\begin{aligned}
y_{n+1} &= y_n + \frac32 \Delta t k_1 + \frac12 \Delta t k_2, \\
(I-\gamma\Delta t\widehat{A})k_1 &= f(t_n,y_n), \\
(I-\gamma\Delta t\widehat{A})k_2 &= f(t_{n+1},y_n + \Delta t k_1) - 2k_1.
\end{aligned}
\end{equation}
% where $\widehat{A}\approx A$.
The method is second order consistent for any $\widehat{A}\in\Rrnn$
and, to have good stability properties, one usually takes $\widehat{A}\approx A$.
Typically, $\widehat{A}$ corresponds to the terms in $A$ which have
to be integrated implicitly in time.  For instance in~\cite{ROS2}, for
advection-diffusion-reaction problems,
$\widehat{A}$ is taken such that
$$
I-\gamma\Delta t\widehat{A}=
(I-\gamma\Delta t \Adiff)(I-\gamma\Delta t A_{\mathrm{react}}),
$$
where $A_{\mathrm{diff}}$ contains diffusion terms and $A_{\mathrm{react}}$
is the reaction Jacobian.
In this work we take $\widehat{A}$ to be either $A$ or the diffusion part of $A$.
Following suggestion in~\cite[Chapter~IV.5, Remark 5.2]{HundsdorferVerwer:book}
we set $\gamma=1$.

\section{Numerical experiments}
\label{s:num}
Numerical experiments described here are carried in Matlab
on a Linux PC with 6 Intel Core i5-8400 CPUs of 2.80GHz,
with 16 Gb RAM.

%---------------------------------------------------------------------
\subsection{Test 1: time dependent source and boundary conditions}
\begin{figure}[t]%
\centering{%
\includegraphics[width=0.48\linewidth]{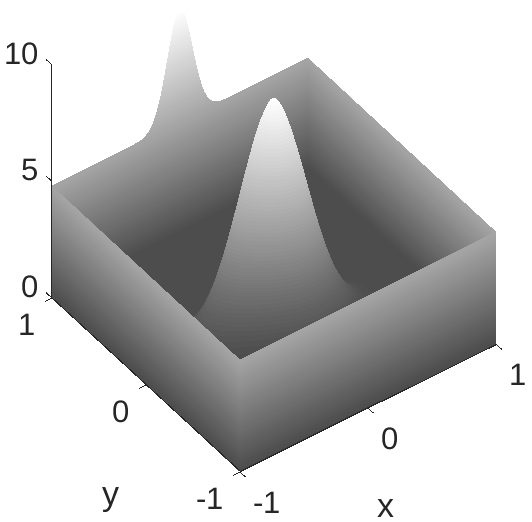}%
\includegraphics[width=0.48\linewidth]{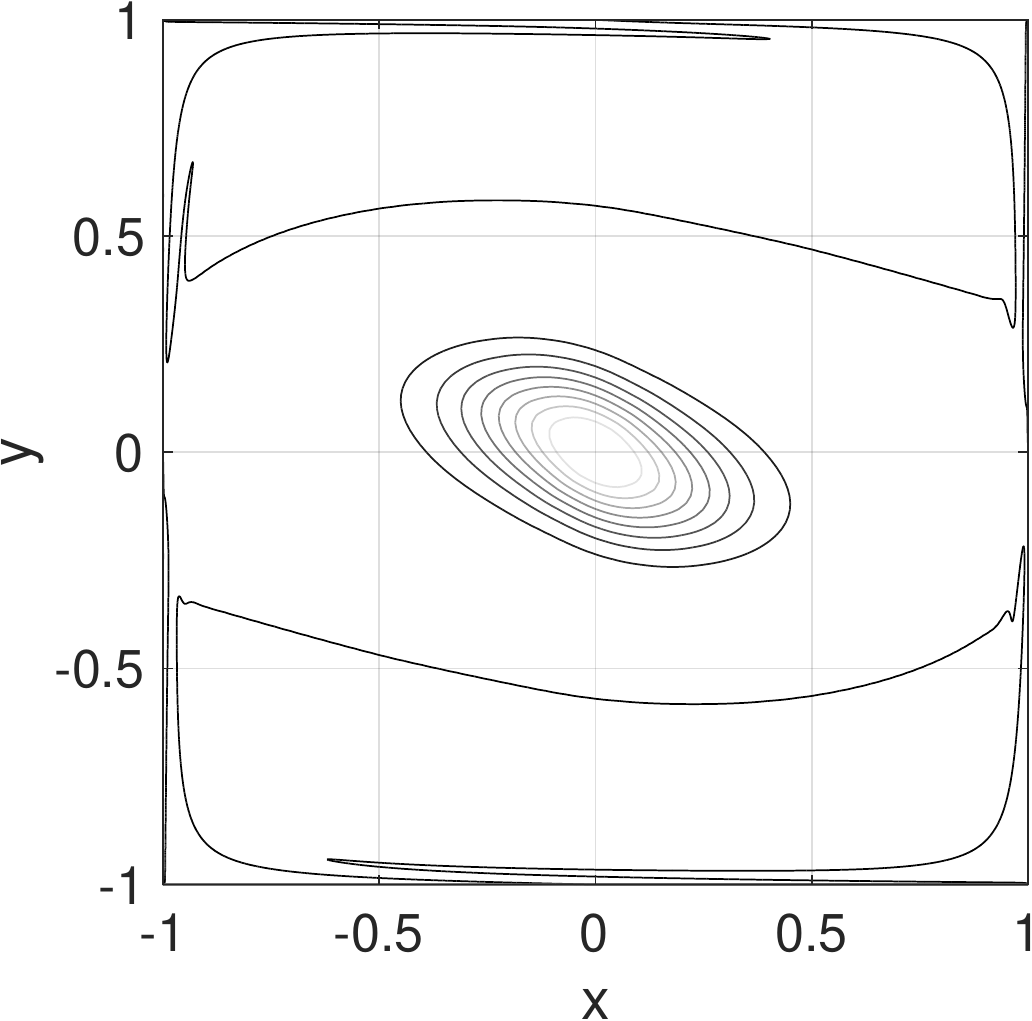}}
\caption{Solution function~\eqref{yex1} on the mesh $256\times 256$
at final time $T=1000$ as surface (left) and contour (right) plots.}
\label{f:yex1}
\end{figure}

In this test we solve~\eqref{ivp} where $A$ is a finite-element discretization
of the two-dimensional convection--diffusion problem:
\begin{equation}
\label{bp_ifiss}  
  -\nu\nabla^2u  + \bm{v}\cdot\nabla u = 0,\qquad
  u=u(x,y),\quad (x,y)\in[-1,1]\times[-1,1],
\end{equation}
where $\nu$ is the viscosity parameter and
the velocity field is $\bm{v}=[ v_1(x,y), v_2(x,y) ]$, 
$$
v_1(x,y) = y(1-x^2), \qquad
v_2(x,y) = x(y^2-1),
$$
For this test, the function $g(t)$ in~\eqref{ivp} takes the form
$$
g(t)\equiv \yex'(t) + A \yex(t),
$$
where $\yex(t)$ is exact solution function chosen as
\begin{equation}
\label{yex1}
\begin{aligned}
\yex(t)   &= \alpha(t)\bigl( A^{-1}\gbc + T\varphi(-TA)\gpeak \bigr),
\\
\alpha(t) &= 1 - e^{-t/300} + e^{-t/100}.
\end{aligned}  
\end{equation}
Here $\gbc\in\Rr^N$ is a vector containing Dirichlet boundary values
prescribed below
and the vector $\gpeak\in\Rr^N$ consists of the values of function
$e^{-10x^2 -50y^2}$ on the mesh.
The boundary conditions imposing by $\gbc$ are 
$$
u(-1,y)=u(1,y)=u(x,-1)=5,
\qquad u(x,1)=5 + 5e^{-50x^2}.
$$
Note that $A^{-1}\gbc$ is the steady state solution of~\eqref{ivp} for
$g(t)\equiv\gbc$ and $T\varphi(-TA)\gpeak$ is the solution
of~\eqref{ivp} with $g(t)\equiv\gpeak$ at time $t=T$.
The final time is $T=1000$ in this test.  In Figure~\ref{f:yex1}
the exact solution $\yex(T)$ is plotted.

In this test the IFISS finite element
discretization~\cite{IFISS,IFISS_sirev} by
bilinear quadrilateral ($Q_1$) finite elements with 
the streamline upwind Petrov--Galerkin (SUPG) stabilization
is employed.
We set viscosity to $\nu=1/6400$ and use nonuniform Cartesian stretched
$256\times 256$ and $512\times 512$ grids with default refinement parameters,
which get finer near the domain boundaries, see Table~\ref{t:mesh}.

\begin{table}
\caption{Parameters of the IFISS stretched meshes}
\label{t:mesh}  
\centering{\begin{tabular}{c@{\extracolsep{1em}}cccc}
\hline
mesh   &  $\min h_x=\min h_y$   & $\max h_x=\max h_y$ & ratio  & max.elem.\ Peclet\\
&                        &                    & $\max h_{x,y}/\min h_{x,y}$

& for $\nu=1/6400$\\ 
\svhline
$256\times 256$  & \num{5.9804e-04} & \num{0.0312} & \num{52.17} & \num{1.9989e+02}
\\
$512\times 512$  & \num{2.0102e-04} & \num{0.0176} & \num{87.5535} & \num{1.1248e+02}
\\\hline
\end{tabular}}
\end{table}

When constructing an advection-diffusion
matrix, the IFISS package provides the value the maximum finite element
grid Peclet number, which is evaluated as
$$
\frac{1}{2\nu}\min\left\{\frac{h_x}{\cos\alpha},
\frac{h_y}{\sin\alpha}\right\}\|\bm{v}\|_2,
\quad \alpha=\arctan\frac{v_2}{v_1},
$$
is $\approx 25$, where $h_{x,y}$ and $v_{1,2}$ are respectively
the element sizes and the velocity components.
The maximum element Peclet numbers reported for these meshes
are given in Table~\ref{t:mesh}.
Due to the SUPG stabilization, the resulting matrices for both meshes
are weakly nonsymmetric: the ratio $\|A-A^T\|_1/\|A+A^T\|_1$
amounts approximately to $0.022$ (mesh $256\times 256$)
and $0.012$ (mesh $512\times 512$).

In addition to the requested accuracy tolerance, two input parameters have
to provided to EBK: the number of the truncated SVD terms $m$ and the number
of time snapshots $n_s$ to construct approximation~\eqref{g=Up}.
From the problem description, we see that $\yex(t)$ is a linear
combination of two linearly independent vectors for any $t$.
Hence, $g(t)$ is a linear combination of no more than four vectors
and we should take $m\leqs 4$.  The actual situation is displayed
by the singular values available from the thin SVD of the time samples:
the largest truncated singular value $\sigma_{m+1}$ is an upper bound
for the truncation error $\|g(t)-Up(t)\|$, see, e.g.,~\cite{GolVanL}.
In this case it turns out that taking $m=2$ is sufficient.
A proper snapshot number $n_s$ can be estimated from given $\alpha(t)$
or by checking, for constructed $U$ and $p(t)$, the actual error
$\|g(t)-Up(t)\|$ a posteriori, see Table~\ref{t:Upt}.
Based on this, we set $n_s=120$ in all EBK runs in this test.
This selection procedure for $n_s$ is computationally very
cheap and can be done once, before all the test runs.

\begin{table}[t]
  \caption{Error of approximation~\eqref{g=Up}, $256\times 256$ mesh.
  The EBK errors are obtained for $\tol=10^{-6}$.}\label{t:Upt}
% In test_expm_krylov_block_g2.m set
% options.check_svd = 1; 
% The run
% test_expm_krylov_block_g2(1e-6,100,30)
% test_expm_krylov_block_g2(1e-6,100,60)
% test_expm_krylov_block_g2(1e-6,100,120)
\centering{\begin{tabular}{c@{\extracolsep{1em}}ccc}
\hline
$n_s$   & $\displaystyle\max_{s\in[0,T]}\|g(s)-Up(s)\|$
                        & $\dfrac{\int_0^T\|g(s)-Up(s)\|ds}
    {\int_0^T\|g(s)\|ds}$   &
    EBK error~\eqref{err1}\\
    \svhline
    30  & {\tt2.37e-03} & {\tt1.82e-05} & {\tt2.24e-05}  \\
    60  & {\tt1.27e-04} & {\tt9.83e-07} & {\tt1.23e-06}  \\
    120 & {\tt7.40e-06} & {\tt5.73e-08} & {\tt7.90e-08}  \\
    \hline
\end{tabular}}
\end{table}

As the problem is two-dimensional, linear systems with $A$
can be solved efficiently by sparse direct methods.
Therefore to solve the linear systems in ROS2, we use the Matlab
standard sparse LU factorization (provided by UMFPACK)
computing it once and using at each time step.

The error reported below for the test runs is measured as
\begin{equation}
  \label{err1}
  \mathrm{error} = \dfrac{\|y(T)-\yex(T)\|}{\|\yex(T)\|}.
\end{equation}
The results of the test runs are presented in Tables~\ref{t:256}
and~\ref{t:512}.
As we see, EBK turns out to be more efficient than the other solvers.
Within the EE2 integrator, the change of the Krylov subspace solver
from EXPOKIT's \texttt{phiv} to the RT-restarted algorithm
leads to a significant increase in efficiency.
Note that this gain is not due to the restarting
but due a more reliable residual-based error control.
Restarting is usually not done because, due to a sufficiently
small $\Delta t$, just a couple Krylov
steps are carried out in EE2 each time step.
In both EE2/EXPOKIT and EE2/RT we should be careful with setting
a proper tolerance value, which is used at each time step for
stopping the Krylov subspace method evaluating 
the $\varphi$ matrix function.  Taking a large tolerance value
may lead to an accuracy loss. For increasingly small
tolerance values the same accuracy will be observed (as it is determined
by the time step size) at a higher cost: more matrix-vector multiplications
per time step will be needed for the $\varphi$ matrix function evaluations.

From Tables~\ref{t:256} and~\ref{t:512} we also see that
the ROS2 solver becomes less efficient than EE2/RT on the finer mesh
as the costs for solving linear systems become
more pronounced.

\begin{table}[t]
\caption{Test~1. Results for the $256\times 256$ mesh}\label{t:256}
% 
% test_ros2a(5,'full') with fname = 'square_cd_0256_nu1_6400.mat';
%
\begin{tabular}{p{4.89cm}p{2cm}p{2.4cm}p{2cm}}
\hline
method & CPU     & fevals$^a$, & error  \\
       & time, s & l.s.s.$^b$  &        \\\svhline
EBK, $\tol=10^{-4}$ %, $n_s=120$
       & 0.36    & 20, ---     & {\tt8.01e-08} \\
EBK, $\tol=10^{-6}$ %, $n_s=120$
       & 0.40    & 24, ---     & {\tt7.90e-08} \\
EE2/RT, $\Delta t = 20$, $\tol=10^{-4}$
       & 1.84    & 500, ---   & {\tt1.51e-03} \\
EE2/RT, $\Delta t = 10$, $\tol=10^{-4}$
       & 3.52    & 900, ---   & {\tt3.79e-04} \\
EE2/RT, $\Delta t = 5$, $\tol=10^{-4}$
       & 6.93    & 1800, ---   & {\tt9.50e-05} \\
EE2/EXPOKIT, $\Delta t = 20$, $\tol=10^{-4}$
       & 17.99    & 9408, ---   & {\tt1.51e-03} \\
EE2/EXPOKIT, $\Delta t = 10$, $\tol=10^{-4}$
       & 24.53    & 12608, ---   & {\tt3.79e-04} \\
EE2/EXPOKIT, $\Delta t = 5$, $\tol=10^{-4}$
       & 37.74    & 19200, ---   & {\tt9.50e-05} \\
ROS2, $\widehat{A}=A$, $\Delta t = 20$ 
       & 1.95    & 100, 100    & {\tt3.03e-03} \\
ROS2, $\widehat{A}=A$, $\Delta t = 10$ 
       & 3.59    & 200, 200    & {\tt7.60e-04} \\
ROS2, $\widehat{A}=A$, $\Delta t = 5$ 
       & 6.85    & 400, 400    & {\tt1.91e-04} \\  
ROS2, $\widehat{A}=\Adiff$, $\Delta t = 2$ 
       & 19.55   & 1000, 1000  & {\tt8.49e-04} \\
ROS2, $\widehat{A}=\Adiff$, $\Delta t = 1$ 
       & 34.72   & 2000, 2000  & {\tt7.59e-06} \\
ROS2, $\widehat{A}=\Adiff$, $\Delta t = 0.5$ 
       & 68.33   & 4000, 4000  & {\tt1.90e-06} \\  
\hline
\end{tabular}
$^a$ number of function evaluations or matvec (matrix-vector) products\\
$^b$ number of linear system solutions
\end{table}

\begin{table}[t]
\caption{Test~1. Results for the $512\times 512$ mesh}\label{t:512}
% 
% test_ros2a(5,'full') with fname = 'square_cd_0512_nu1_6400.mat';
%
\begin{tabular}{p{4.89cm}p{2cm}p{2.4cm}p{2cm}}
\hline
method & CPU     & fevals$^a$, & error  \\
       & time, s & l.s.s.$^b$  &        \\\svhline
EBK, $\tol=10^{-4}$
       & 1.26    & 4, ---     & {\tt3.08e-08} \\
EBK, $\tol=10^{-6}$
       & 1.31    & 8, ---     & {\tt2.33e-08} \\
EE2/RT, $\Delta t = 20$, $\tol=10^{-4}$
       &  9.10   & 450, ---    & {\tt8.91e-04} \\
EE2/RT, $\Delta t = 10$, $\tol=10^{-4}$
       & 17.97   & 900, ---    & {\tt2.40e-04} \\
% test_expm_krylov_block_g2(1e-4,200,120)
EE2/RT, $\Delta t = 5$, $\tol=10^{-4}$
       & 35.90   & 1800, ---   & {\tt8.07e-05} \\
ROS2, $\widehat{A}=A$, $\Delta t = 20$ 
       & 11.82   & 100, 100    & {\tt1.90e-03} \\
ROS2, $\widehat{A}=A$, $\Delta t = 10$ 
       & 22.68   & 200, 200    & {\tt5.46e-04} \\
ROS2, $\widehat{A}=A$, $\Delta t = 5$ 
       & 36.91   & 400, 400    & {\tt1.85e-04} \\
ROS2, $\widehat{A}=\Adiff$, $\Delta t = 2$ 
       & 86.55   & 1000, 1000  & {\tt5.73e-03} \\
ROS2, $\widehat{A}=\Adiff$, $\Delta t = 1$ 
       & 167.95  & 2000, 2000  & {\tt4.44e-06} \\
ROS2, $\widehat{A}=\Adiff$, $\Delta t = 0.5$ 
       & 331.38  & 4000, 4000  & {\tt1.11e-06} \\  
\hline
\end{tabular}
$^a$ number of function evaluations or matvec products\\
$^b$ number of linear system solutions
\end{table}

%---------------------------------------------------------------------
\subsection{Test 2: time dependent boundary conditions}
In the previous test we see that the EBK solver apparently profits from
the specific source function, exhibiting a very quick convergence.
Although this is not an unusual situation, we now consider another
test problem which appears more difficult for EBK.
We take the same matrix~$A$ as in the first test and the following
initial value vector $v$ and source function $g(t)$:
$$
%\begin{aligned}
%\end{aligned}
g(t)=\alpha(t)\gbc, \qquad
v = -T\varphi(-TA)\gpeak,
$$
where $\alpha(t)$ and $\gbc$ are the same as in~\eqref{yex1}.
This test problem does not have a known analytical solution
and we compute a reference solution $\yref(t)$ by running EE2/RT with
a tiny time step size.  The errors of computed numerical solutions
$y(t)$ reported below are
$$
\textrm{error} = \frac{\|y(T)-\yref(T)\|}{\|\yref(T)\|}
$$
Note that $\yref(t)$ is influenced by the same space error as $y(t)$,
hence, the error shows solely the time error.

From the problem definition we see that the number of SVD terms $m$
can be at most $2$.  Therefore, in this test 
EBK is run with the block size $m=2$ and
$n_s=80$ time snapshots (the value is determined in the same
way as in Test~1).
For this test we include in comparisons the two solvers which
come out as best in the first test, EBK and EE2/RT. 
The results presented in Table~\ref{t:256.2} show
that EBK does require more steps for this test
but is still significantly more efficient than EE2/RT.

\begin{table}[t]
\caption{Test 2.  Results for the $256\times 256$ mesh}\label{t:256.2}
% 
% test_ros2a(5,'full') with fname = 'square_cd_0512_nu1_6400.mat';
%
\begin{tabular}{p{4.89cm}p{2cm}p{2.4cm}p{2cm}}
\hline
% test_expm_krylov_block_g3(1e-8,100,120,2)  
method & CPU     & fevals$^a$, & error  \\
       & time, s & l.s.s.$^b$  &        \\\svhline
EBK, $\tol=10^{-4}$, $n_S=80$
       & 0.77    & 36, ---     & {\tt1.83e-05} \\
EBK, $\tol=10^{-6}$, $n_S=80$
       & 1.32    & 50, ---     & {\tt1.91e-07} \\
EE2/RT, $\Delta t = 10$, $\tol=10^{-6}$
       & 6.53    & 1306, ---   & {\tt8.91e-05} \\
EE2/RT, $\Delta t = 5$, $\tol=10^{-6}$
       & 11.54   & 2406, ---   & {\tt5.58e-05} \\
\hline
\end{tabular}
$^a$ number of function evaluations or matvec products\\
$^b$ number of linear system solutions
\end{table}

\section{Conclusions}
\label{s:concl}
We show that exponential time integrators can be an attractive
option for integrating advection--diffusion problems in time,
as they possess good accuracy as well as stability properties.
In presented tests, they outperform state-of-the-art implicit-explicit
ROS2 solvers. 
Exponential solvers which are able to exploit their matrix function
evaluation machinery for a whole time interval (such as EBK in
this paper) appear to be preferable to exponential integrators where
matrix functions have to be evaluated at each time step.

%% \begin{figure}[t]%
%% \sidecaption[t]%
%% % Use the relevant command for your figure-insertion program
%% % to insert the figure file.
%% % For example, with the option graphics use
%% \includegraphics[scale=.65]{figure}%
%% %
%% % If no graphics program available, insert a blank space i.e. use
%% %\picplace{5cm}{2cm} % Give the correct figure height and width in cm
%% %
%% %\caption{Please write your figure caption here}
%% \caption{If the width of the figure is less than \SI{7.8}{\cm} cm
%% use the \texttt{sidecapion} command to flush the caption on the
%% left side of the page. If the figure is positioned at the top of
%% the page, align the sidecaption with the top of the figure -- to
%% achieve this you simply need to use the optional argument \texttt{[t]}
%% with the \texttt{sidecaption} command}\label{fig:3}% Give a unique label
%% \end{figure}
%% % or

%---------------------------------------------------------------------
%\begin{acknowledgement}
%\end{acknowledgement}

\bibliographystyle{spmpsci}
\bibliography{matfun,my_bib,split}

\end{document}